\RequirePackage{fix-cm}
\documentclass[smallextended]{svjour3}       
\smartqed  
\usepackage{graphicx}
\usepackage{amsmath}
\usepackage{amssymb}
\usepackage{bm}
\usepackage{float}
\usepackage{multirow}

\DeclareMathOperator{\E}{E} \DeclareMathOperator{\var}{var}
\DeclareMathOperator{\bias}{bias} \DeclareMathOperator{\MSE}{MSE}
\newcommand{\bR}{\mathbb{R}}

\newenvironment{problist}
   {
      \begin{list}
         {---}
         {
            \setlength{\itemsep}{.5ex}
            \setlength{\parsep}{0ex}
            \setlength{\parskip}{0ex}
            \setlength{\topsep}{.5ex}
         }
   }
   {
      \end{list}
   }
\begin{document}

\title{Reduced bias nonparametric lifetime density and hazard estimation
}


\author{Arthur Berg \and Dimitris Politis \and Kagba Suaray \and Hui Zeng
}


\institute{Arthur Berg \at
              Penn State College of Medicine \\
              Division of Biostatistics \& Bioinformatics\\
              \email{berg@psu.edu}           
           \and
           Dimitris Politis \at
              University of California, San Diego\\
              Department of Mathematics
              \and
              Kagba Suaray \at
              California State University Long Beach\\
              Department of Mathematics \& Statistics
              \and
              Hui Zeng \at 
              Penn State College of Medicine\\
              Division of Biostatistics \& Bioinformatics
}

\date{Received: date / Accepted: date}

\maketitle

\begin{abstract}
Kernel-based nonparametric hazard rate estimation is considered with a special class of infinite-order kernels that achieves favorable bias and mean square error properties.  A fully automatic and adaptive implementation of a density and hazard rate estimator is proposed for randomly right censored data.  Careful selection of the bandwidth in the proposed estimators yields estimates that are more efficient in terms of overall mean squared error performance, and in some cases achieves a nearly parametric convergence rate.  Additionally, rapidly converging bandwidth estimates are presented for use in second-order kernels to supplement such kernel-based methods in hazard rate estimation. Simulations illustrate the improved accuracy of the proposed estimator against other nonparametric estimators of the density and hazard function.  A real data application is also presented on survival data from 13,166 breast carcinoma patients.
\end{abstract}

\section{Introduction}

Hazard rate estimation has been extensively studied in the literature as it encompasses fundamental characteristics of time-to-event data with applications spanning medicine, engineering, and economics.  The first kernel-based nonparametric estimator of the hazard function with non-censored data appeared in \cite{Watson:1964aa}.  For censored data, density estimation approaches are described in \cite{Foldes:1981aa} and \cite{Padgett1984}, and an empirical hazard approach is described in \cite{Tanner:1983aa}.  Kernel-based estimation of the hazard function under censoring was studied by \cite{Yandell1983}, \cite{Ramlau-Hansen1983}, and \cite{Tanner1984}, \cite{Muller:1994aa}, among others.  However, all of these kernel-based approaches capitalized on traditional theory of second-order kernels when constructing their kernel-based estimates.  Through the use of infinite-order kernels, we demonstrate that considerable asymptotic improvements are attainable.   

The benefit of using infinite-order kernels, also called superkernels, in estimating the probability density function under iid data is well known; cf.\ \cite{Devroye92}. More recently, Politis and others have investigated a class of infinite-order kernels that are constructed by taking the Fourier transform of flat-top functions---functions that are flat in a neighborhood of the origin \cite{Berg:2009aa,Berg:2009ab,McMurry:2004aa,Politis99,politis93}. These estimators, under a correctly specified bandwidth, attain mean squared error (MSE) properties superior to their second order analogs and also perform well in small sample simulation studies. These same properties translate nicely to the context of density estimation under random right censoring, as investigated here. Improved MSE convergence rates in nonparametric estimation of the hazard function and derivatives of the density follow as corollaries to the density estimation theory.

In the next section, we define the general class of flat-top infinite-order kernels and, through Theorem \ref{thm: cen-bias}, describe how using
these kernels can cause the bias of density estimators from censored data to become essentially negligible in certain situations.  Section 3 completes the proposed estimator by providing a bandwidth selection algorithm that automatically adapts to the unknown density at hand. A second use of the infinite-order estimators is realized in Section 4 by providing rapidly converging bandwidths for use in second-order kernels. In Section 5, we give practical suggestions for implementing the proposed estimator and provide simulations exhibiting improved performance in estimating the lifetime density and hazard function when compared with other nonparametric estimators including the \texttt{muhaz} and \texttt{pehaz} estimators of \cite{Hess:2014aa,Muller:1994aa} and the \texttt{presmooth} estimator of \cite{Lopez-de-Ullibarri:2013aa}.  In Section 6, the proposed hazard function estimator and the previously mentioned estimators are simultaneously compared on breast carcinoma survival data involving 13,166 women.

\section{Estimation with Infinite-Order Kernels}

We lay out the notation under the context of random right censorship (this can be generalized to allow for left truncation; see for example \cite{cao99}).  Let $X^0_1,\ldots, X^0_n$ be iid lifetime variables with density $f$ and cdf $F$, and independently, let $U_1,\ldots,U_n$ be iid censoring variables with density $g$ and cdf $G$.   We observe the data $Z_i$ and $\Delta_i$ where
\[
Z_i=\min\{X_i^0,U_i\}\quad\text{and}\quad \Delta_i=1_{[X_i^0\le U_i]}\in\{0,1\}
\]
for $i=1,\ldots,n$ (here $1_{[\cdot]}$ represents the indicator function).  We order the pairs $(Z_i,\Delta_i)$ according to the $Z_i$'s and relabel them as $(X_i,\delta_i)$ where $X_i=Z_{(i)}$, the $i^\text{th}$ order statistics of the $Z$'s, and $\delta_i$ is
the indicator variable that accompanies $X_i$, i.e. the concomitant of $X_i$. The Kaplain-Meier estimator is the nonparametric maximum likelihood estimate of the survival function $S(t)=1-F(t)$ given by
\[
\hat{S}(t)=
\begin{cases}
  1, & 0\le t\le X_1\\
  \displaystyle{\prod_{j=1}^{k-1}}\left(\frac{n-j}{n-j+1}\right)^{\delta_j}, &X_{k-1}<t\le X_{k},\ k=2,\ldots,n \\
  0, & t>X_n
\end{cases}
\]
where the height of the jump of $\hat{S}$ at $X_j$ is
\[
s_j=
\begin{cases}
  \hat{S}(X_j)-\hat{S}(X_{j+1}), & j=1,\ldots,n-1\\
  \hat{S}(X_n), & j=n.
\end{cases}
\]
The kernel estimate of $f$ is constructed through the convolution of
$\hat{F}=1-\hat{S}$ with a smooth kernel $K$, i.e.
\begin{equation}
\label{eq: kernel est}
\begin{split}
  \hat{f}(x)&=\frac{1}{h}\int_{-\infty}^\infty K\left(\frac{x-t}{h}\right)\,d\hat{F}(t)=\frac{1}{h}\sum_{j=1}^{n}s_j
  K\left(\frac{x-X_j}{h}\right).
\end{split}
\end{equation}
See \cite{Foldes:1981aa} and \cite{Marron87} for background and properties of this estimator. Many authors require $K$ to be of compact support for ease of
analysis, but this is unnecessary; see for example \cite{Tanner:1983aa}.
Therefore we only assume $K$ is an even function that integrates to
one.

It will be assumed that sufficient conditions are satisfied so that
\begin{equation}
\label{eq: var assumption}
\var\left(\hat{f}(x)\right)=O\left(\frac{1}{nh}\right).
\end{equation}
This typically requires $h\rightarrow0$ as $nh\rightarrow\infty$ and $n\rightarrow\infty$, the lifetime density $f$ to be continuously differentiable at $x$, and the censored distribution to be of compact support.  Under these conditions, precise variance expressions are provided in \cite{Gijbels93}.

Following \cite{politis01}, we now describe a class of
infinite-order kernels constructed from the Fourier transform of a
flat-top function. We start in the Fourier domain with a function
$\kappa$ given by
\begin{equation}
\label{eq: kappa} \kappa(t)=
\begin{cases}
  1, & |t|\le c\\
  g(|t|), & \text{otherwise}
\end{cases}
\end{equation}
where $c$ is any positive constant, and $g$ is any continuous, square-integrable function that is
bounded in absolute value by one and satisfies $g(|c|)=1$. Then
the infinite-order kernel corresponding to $\kappa$ is the Fourier
transform of $\kappa$, specifically,
\begin{equation}
\label{eq: flat-top kernel} K(x)=\frac{1}{2\pi}\int_{-\infty}^\infty
\kappa(t) e^{-i t x} \, dt,
\end{equation}
or equivalently,
\begin{equation}
\label{eq: identity}
\frac{1}{h}K(x/h)=\frac{1}{2\pi}\int_{-\infty}^\infty \kappa(t
h)e^{-itx}\,dt.
\end{equation}

Let $\phi(t)$ be the characteristic function corresponding to
$f(x)$, i.e. $\phi(t)$ is the inverse Fourier transform of $f(x)$
given by
\[
\phi(t)=\int_{-\infty}^\infty e^{i t x}f(x)\,dx.
\]
A natural estimator of the characteristic function is
\begin{equation}
\label{eq:phihat}
\begin{split}
\hat{\phi}(t)&=\int_{-\infty}^\infty e^{i t
x}\,d\hat{F}(x)=\sum_{j=1}^n s_j\, e^{i t X_j}
\end{split}
\end{equation}
In the context of non-censored data, $\hat{\phi}$ is an unbiased estimator of $\phi$, but in the presence of censoring, bias is present.  We assume the bias of $\hat\phi(x)$ is $O(\frac{1}{n})$, which is justified by the following lemma.  
\begin{lemma}
\label{lem:bias}
Suppose $g(x)$ (the censored density) is compactly supported and contains the support of $f(x)$, then
\begin{equation}
\label{eq:phibias}
\text{bias}(\hat{\phi}(x))=O\left(\frac{1}{n}\right).
\end{equation}
\end{lemma}
The assumption in Lemma \ref{lem:bias} that the support of the censored distribution contains the support of the lifetime distribution can be found in other papers (see e.g. \cite{Una-Alvarez:2015aa}).  This assumption may be limiting in certain applications, but often it suffices to estimate a truncated lifetime distribution.

The bias of $\hat f(x)$ is smaller when $f(x)$ is smoother (has more derivatives), and the more smooth $f(x)$ is, the faster its characteristic function decays to zero.   The following assumptions classify the smoothness of $f(x)$ into one of three categories.
\begin{itemize}
  \item\textbf{Assumption $\bm{A(r)}$:}  There is an $r>0$ such that
  $\int_{-\infty}^\infty |t|^r\,|\phi(t)|\,dt <\infty$. 
  \item\textbf{Assumption $\bm{B}$:}  There are positive constants $d$ and $D$ such
  that $|\phi(t)|\le D e^{-d|t|}$.
  \item\textbf{Assumption $\bm{C}$:}  There is a positive constant $b$ such that
$|\phi(t)|=0$ when $|t|\ge b$.
\end{itemize}

The following theorem provides rates of the bias and MSE under each of the assumptions above.
\begin{theorem}
\label{thm: cen-bias} Suppose $\hat{f}(x)$ is the kernel estimator
as defined in \eqref{eq: kernel est} with infinite-order kernel
given by \eqref{eq: flat-top kernel} and assume the variance
assumption in \eqref{eq: var assumption} and the bias assumption in \eqref{eq:phibias}.  
\begin{problist}
  \item[(i)]  Suppose assumption $A(r)$ holds.  Let $h\sim a
  n^{-\beta}$ (for any $a>0$) with $\beta=(2r+1)^{-1}$, then
\[
\sup_{x\in\bR}
\left|\bias\left\{\hat{f}(x)\right\}\right|=o\left(n^{\frac{-r}{2r+1}}\right)\quad\text{and}\quad
\MSE\left\{\hat{f}(x)\right\}=O\left(n^{\frac{-2r}{2r+1}}\right).
\]
  \item[(ii)]  Suppose assumption $B$ holds.  Let $h\sim 1/(a\log n)$ with $a> 1/(2d)$, then
\[
\sup_{x\in\bR}
\left|\bias\left\{\hat{f}(x)\right\}\right|=O\left(\frac{1}{\sqrt{n}}\right)\quad\text{and}\quad
\MSE\left\{\hat{f}(x)\right\}=O\left(\frac{\log n}{n}\right).
\]
\item[(iii)]  Suppose assumption $C$ holds.  Let $h\le 1/b$, then
\[
\sup_{x\in\bR}
\left|\bias\left\{\hat{f}(x)\right\}\right|=O\left(\frac{1}{n}\right)\quad\text{and}\quad
\MSE\left\{\hat{f}(x)\right\}=O\left(\frac{1}{n}\right).
\]
\end{problist}
\end{theorem}

This theorem illustrates the mean square error of $\hat{f}(x)$ is just as good as second-order kernel density estimation when $f(x)$ is only twice differentiable ($r=2$), but considerable improvements are gained when more smoothness of $f(x)$ is present.  Even a parametric convergence rate is possible when assumption $C$ is satisfied.  Parametric convergence rates in non-censored data has also demonstrated by others \cite{Davis:1977aa,Ibragimov:1983aa,Chacoon:2007aa}, such as for data following a Vall\'{e}-Poussin density given by
\[
f(x)=\frac{1-\cos x}{\pi x^{2}},\quad x\in \mathbb{R}
\]
with finitely-supported characteristic function given by
\[
\phi(t)=(1-|t|)I_{[-1,1]}(t).
\]

\begin{corollary}  The hazard function $\lambda(x)=f(x)/S(x)$ is naturally
estimated by $\hat{\lambda}(x)=\hat{f}(x)/\hat{S}(x)$, and since $\hat{S}$
is a $\sqrt{n}$-convergent estimator of $S$, this estimate of the
hazard function has the same MSE convergence rates as $\hat{f}$ in
the above theorem.  Specifically:
\begin{problist}
  \item[(i)]  Under assumption $A(r)$,
$\MSE(\hat{\lambda}(x))=O\left(n^{\frac{-2r}{2r+1}}\right)$;
  \item[(ii)]  Under assumption $B$,
$\MSE(\hat{\lambda}(x))=O\left(\frac{\log n}{n}\right)$;
  \item[(iii)]  Under assumption $C$,
$ \MSE(\hat{\lambda}(x))=O\left(\frac{1}{n}\right)$.
\end{problist}
\end{corollary}

Additionally, the $p^{\text{th}}$ derivative of $f$ can be estimated by the
$p^{\text{th}}$ derivative of $\hat{f}(x)$; i.e. if $K^{(p)}(x)$ is
the $p^{\text{th}}$ derivative of $K(x)$, then
\begin{equation}
\label{eq: fhatp} \hat{f}_p(x)=\frac{1}{h^{p+1}}\sum_{j=1}^{n}s_j
  K^{(p)}\left(\frac{x-X_j}{h}\right)
\end{equation}
is an estimate of the $p^{\text{th}}$ derivative of $f$ \cite{Singh:1977aa}.  Similarly, under sufficient conditions on $f$, the variance of this
estimator is
\begin{equation}
\label{eq: var2}
\var\left(\hat{f}_p(x)\right)=O\left(\frac{1}{n\,h^{p+1}}\right).
\end{equation}
The previous theorem is now generalized in the following theorem
to give asymptotic bias and MSE rates of $\hat{f}_p(x)$ with
infinite-order kernels.

\begin{theorem}
\label{thm: derivative} \label{thm: cen-bias2} Suppose
$\hat{f}_p(x)$ is the kernel estimator as defined in \eqref{eq:
fhatp} where $K$ is an infinite-order kernel, and assume \eqref{eq:phibias} and \eqref{eq: var2} hold.
\begin{problist}
  \item[(i)]  Suppose assumption $A(r+p)$ holds.  Let $h\sim a
  n^{-\beta}$ (for any $a>0$) with $\beta=(2r+p+1)^{-1}$, then
\[
\sup_{x\in\bR}
\left|\bias\left\{\hat{f}_p(x)\right\}\right|=o\left(n^{\frac{-r}{2r+p+1}}\right)\quad\text{and}\quad
\MSE\left\{\hat{f}_p(x)\right\}=O\left(n^{\frac{-2r}{2r+p+1}}\right).
\]
  \item[(ii)]  Suppose assumption $B$ holds.  Let $h\sim 1/(a\log n)$ with $a> 1/(2d)$, then
\[
\sup_{x\in\bR}
\left|\bias\left\{\hat{f}_p(x)\right\}\right|=O\left(\frac{1}{\sqrt{n}}\right)\quad\text{and}\quad
\MSE\left\{\hat{f}_p(x)\right\}=O\left(\frac{\log n}{n}\right).
\]
\item[(iii)]  Suppose assumption $C$ holds.  Let $h\le 1/b$, then
\[
\sup_{x\in\bR}
\left|\bias\left\{\hat{f}_p(x)\right\}\right|=O\left(\frac{1}{n}\right)\quad\text{and}\quad
\MSE\left\{\hat{f}_p(x)\right\}=O\left(\frac{1}{n}\right).
\]
\end{problist}
\end{theorem}

In particular, we see that if the underlying density is infinitely
smooth (as in the case of assumptions $B$ and $C$), then the same
asymptotic MSE rates of $\hat{f}_p(x)$ hold for every $p$.

The bias properties stated in Theorem 2.1 and Theorem 2.3 are consistent with the properties of other kernel-based estimators with infinite-order kernels used in different contexts; see e.g. \cite{Berg:2009ab,Politis99}.  

\section{Bandwidth Selection Algorithm}
\label{sec:abc}
Theorem \ref{thm: cen-bias} in the previous section assumes one is handed a bandwidth $h$ that is precisely molded to the underlying smoothness of the density of interest $f(x)$.  In general, however, one does not necessarily know the level of smoothness of the underlying density.  This section presents a simple algorithm, adapted from \cite{Politis03}, that automatically adjusts to the unknown smoothness of $f(x)$.  This bandwidth estimation procedure is consistent for the optimal bandwidth $h$ under assumptions B and C, and under assumption A($r$), the bandwidth algorithm still adapts to the underlying smoothness but consistency does not hold.   

Let $\hat{\phi}(t)$ be the estimate of the characteristic function as given in \eqref{eq:phihat}.  The algorithm, in essence, searches for the smallest value $t^*$ such that $\hat{\phi}(t)\approx 0$ in which case the bandwidth estimate is taken to be $\hat{h}\approx1/t^*$. The specific details are provided in the following algorithm.

\begin{quote}
  \textsc{Bandwidth Selection Algorithm}\\
  Let $C>0$ be a fixed constant, and $\varepsilon_n$ be a nondecreasing
  sequence of positive real numbers tending to infinity such that $\varepsilon_n=o(\log
  n)$.
  Let $t^*$ be the smallest number such that
  \begin{equation}
  \label{eq: alg}
|\hat{\phi}(t)|<C\sqrt{\frac{\log_{10} n}{n}}\qquad\text{ for all
}t\in (t^*,t^*+\varepsilon_n)
  \end{equation}
  Then let $\hat{h}=c/t^*$ where $c$ is the ``flat-top radius'' given in \eqref{eq: kappa}.
\end{quote}

\begin{remark}
\label{rmk:C} The positive constant $C$ and the choice of sequence $\varepsilon_{n}$ are irrelevant in the asymptotic theory, but certainly relevant for finite-sample calculations.  The main idea behind the algorithm is to determine the smallest $t$ such that $\phi(t)\approx 0$, and in most cases this can be visually seen without explicitly providing the quantities $C$ and $\varepsilon_{n}$ in \eqref{eq: alg}. 
\end{remark}
\begin{remark}
If $g(t)$ in \eqref{eq: kappa} is very close to
one, in a neighborhood of the type $[c,c+ \eta]$, then the
``flat-top radius'' is effectively increased to $c+\eta$.
In this case, we would let $\hat{h}=(c+\eta)/t^*$ in the bandwidth
selection algorithm.  This is particularly relevant when considering infinitely smooth flat-top functions \cite{McMurry:2004aa}.
\end{remark}

\begin{theorem}
\label{thm: cen-band} Assume the following two assumptions on $\hat\phi(t)$:
\begin{equation}
\label{eq: assumption1}
\max_{s\in(0,1)}|\hat{\phi}(t+s)-\phi(t+s)|=O_P(1/\sqrt{n})
\end{equation}
and
\begin{equation}
\label{eq: assumption2}
\max_{s\in(0,n)}|\hat{\phi}(t+s)-\phi(t+s)|=O_P\left(\frac{\log
n}{\sqrt{n}}\right)
\end{equation}
uniformly in $t$.
\begin{itemize}
\item[(i)] If $|\phi(t)|\sim A|t|^{-d}$ for some positive constants $A$ and
$d$, then
\[
\hat{h}\stackrel{P}{\sim} \tilde{A}\left(\frac{\log
n}{n}\right)^{\frac{1}{2d}};
\]
here $A\stackrel{P}{\sim}B$ means $A/B\rightarrow 1$ in probability.
\item[(ii)]  If $|\phi(t)|\sim A \xi^{|t|}$ for some $\xi\in(0,1)$ and $A>0$, then
\[
\hat{h}\stackrel{P}{\sim} 1/(\tilde{A}\log n).
\]
where $\tilde{A}=-1/\log \xi$.
\item[(iii)]  If $|\phi(t)|=0$ when $|t|\ge b$, then $\hat{h}\stackrel{P}{\sim} 1/b$.
\end{itemize}
\end{theorem}
\begin{remark}
The two assumptions \eqref{eq: assumption1} and \eqref{eq: assumption2} are typical assumptions invoked in this type of an algorithm (see e.g. \cite{Politis03}), and verification of these assumptions, particularly with censored data, can be difficult and is not pursued here.    

\end{remark}
Theorem \ref{thm: cen-band} shows that the proposed bandwidth
selection algorithm adapts to the underlying degree of smoothness of
the density, yielding bandwidth estimates that largely match the ideal bandwidths in
Theorem \ref{thm: cen-bias}. When there is only polynomial decay of
the characteristic function, as in part (i) of the above theorem,
the bandwidth selection algorithm produces a slightly smaller
bandwidth than the theoretically optimal bandwidth given in Theorem
\ref{thm: cen-bias}, but the discrepancy diminishes with faster
decay.

\section{Bandwidth Selection for $2^{\text{nd}}$-Order Kernels}  We now propose a bandwidth selection procedure for use with
\textit{second-order} kernels, based on using the infinite-order
estimators as pilots in the plug-in approach to bandwidth selection.
Although Theorem \ref{thm: cen-bias} demonstrates asymptotic superiority of
using infinite-order kernels over second order kernels, the choice of bandwidth in estimation may be more critical than the choice of kernel.  This hybrid approach provides improved (rapidly converging) bandwidth estimates for kernel density estimators using $2^{\text{nd}}$-order kernels.  

We begin with expressions for the MSE and the mean integrated square error (MISE) of
$\hat{f}(x)$ with a symmetric second-order kernel $\Lambda$ and standard
assumptions on $f$ and $G$.  The MISE below is slightly generalized to incorporate a nonnegative weight function $\omega(x)$ to control the influence of error in the tails of the estimated density.  The MSE and MISE calculations will assume the following conditions:
\begin{itemize}
\item[(i)]  $f(x)$, $G(x)$, and $\omega(x)$ are twice differentiable with bounded third derivative in a neighborhood of $x$.
\item[(ii)] $\omega(x)$ is compactly supported whose support is contained inside the support of the censoring distribution.
\item[(iii)]  $\Lambda$ is three times continuously differentiable, its first derivative is integrable, and 
\[
\lim_{|x|\rightarrow\infty}x^{j}\Lambda^{(j)}(x)=0 \qquad(j=0,1,2,3).
\]
The MSE and MISE expressions are now presented; derivations are detailed in \cite{cao99,suaray04}. 
\end{itemize}
\[
\begin{split}
\textrm{MSE}(\hat{f}(x))&=h^4\cdot \left(\frac{f''(x)}{2}c_\Lambda \right)^2 \\
&+\frac{1}{n\,h}\cdot \frac{f(x)}{1-G(x)} d_\Lambda\\
&+\frac{1}{n}\cdot f(x)^2\left[\int_{-\infty}^x \frac{f(r)}{1-G(r)}\,dr-\frac{1}{(1-F(x))(1-G(x))}\right]\\
&+O(h^6)+O\left(\frac{h}{n}\right)+o\left(\frac{1}{n\, h}\right)
\end{split}
\]
where
\[
c_\Lambda=\int_{-\infty}^\infty x^2 \Lambda(x)\,dx \quad and \quad d_\Lambda=\int_{-\infty}^\infty \Lambda^2(x)\,dx
\]
and
\[
MISE(\hat{f})=\int_{-\infty}^\infty MSE(\hat{f}(x))dx.
\]

Minimizing the asymptotically dominant terms in the above expressions with respect to $h$ yields pointwise and globally optimal bandwidths, respectively, given by

\[
h_{\text{MSE}}=\left(\frac{f(x) d_\Lambda}{1-G(x)\left(f''(x) c_\Lambda \right)^2}\right)^{1/5}n^{-1/5}
\]
\[
h_{\text{MISE}}=\left(\frac{d_\Lambda \int_{-\infty}^\infty\frac{f(x)}{1-G(x)}\omega(x)\,dx}{c_\Lambda^2 \int_{-\infty}^\infty
f''(x)^2\omega(x)\,dx}\right)^{1/5}n^{-1/5}
\]

These optimal bandwidths involve values of the unknown values $f(x)$, $f''(x)$ and $G(x)$.  Therefore to estimate the respective bandwidths, we replace these unknown quantities with pilot estimates;  $f(x)$ and $f''(x)$ are replaced with the infinite-order kernel estimates $\hat{f}(x)$ and $\hat{f}_2(x)$
respectively, and $1-G(x)$, the survival function of the \textit{censored} random variables,
is estimated using the Kaplan-Meier estimator with $\Delta_i$ replaced with $1-\Delta_i$.  The bandwidth used in estimating $\hat{f}_2(x)$, and
in general for $\hat{f}_p(x)$, is the same as that derived from
the bandwidth selection algorithm above. If $f(x)$ is sufficiently smooth (for instance when assumption B or C holds), then this bandwidth choice is optimal.   Let $\hat{h}_{\text{MSE}}$  and
$\hat{h}_{\text{MISE}}$ refer to the plug-in estimates corresponding
to $h_{\text{MSE}}$ and $h_{\text{MISE}}$ respectively.  These
estimators have rapid convergence rates due to the ultra-fast
convergence of the plug-in infinite-order kernel estimators, as detailed in the following theorem.
\begin{theorem}
\label{thm: cen-pilot} Assume the conditions of Theorem \ref{thm:
cen-band}, and assume conditions strong enough to ensure \eqref{eq:
var2} holds for $p=2$. Let $\hat{h}_{\text{M}}$ be either
$\hat{h}_{\text{MSE}}$ or $\hat{h}_{\text{MISE}}$ with
$h_{\text{M}}$ being the corresponding $h_{\text{MSE}}$ or
  $h_{\text{MISE}}$.
\begin{itemize}
\item[(i)] If $|\phi(t)|\sim A|t|^{-d}$ for some positive constants $A$ and
$d>3$, then
\[
\hat{h}_{\text{M}}=h_{\text{M}}\left(1+O_p\left(\frac{\log
n}{n}\right)^{\frac{\lceil d-4\rceil}{2d}}\right).
\]
\item[(ii)]  If $|\phi(t)|\sim A \xi^{|t|}$ for some $\xi\in(0,1)$ and $A>0$, then
\[
\hat{h}_{\text{M}}=h_{\text{M}}\left(1+O_p\left(\frac{\log
n}{n}\right)^{\frac{1}{2}}\right).
\]
\item[(iii)]  If $|\phi(t)|=0$ when $|t|\ge b$, then
\[
\hat{h}_{\text{M}}=h_{\text{M}}\left(1+O_p\left(\frac{1}{\sqrt{n}}\right)\right).
\]
\end{itemize}
\end{theorem}

\cite{Marron87} suggest cross-validation as a
means of minimizing the integrated square error (ISE), but the
approach of minimizing ISE was shown in \cite{hall91} to be less
optimal than minimizing the MISE.  In particular, the relative
convergence rates (as in the above theorem) of the cross-validation
approach in \cite{Marron87} are $n^{-1/10}$, regardless of the degree
of smoothness of $f(x)$.  If one uses the plug-in approach that we
have adopted above but with pilots consisting of second-order
kernels, then the relative convergence rates are at best $n^{-2/5}$,
again, regardless of the degree of smoothness of $f(x)$. All of
these rates are considerably slower than the $n^{-1/2}$ rate
afforded by the proposed procedure under a sufficiently smooth
density $f(x)$ (i.e. when $\phi(t)$ has a rapid decay to zero) as
Theorem \ref{thm: cen-pilot} demonstrates.

\section{Simulations}

Many different choices of ``flat-top'' functions \eqref{eq: kappa} can be used to construct an infinite-order kernel, although highly non-smooth shapes like the rectangle, which gives rise to the sinc kernel, should be avoided due to its large and slowly decaying side lobes.  The trapezoidal window, as suggested in \cite{politis95}, can be viewed as smoothening the rectangular window and has more rapidly decaying side lobes.  Another possibility is the infinitely smooth trapezoidal flat-top shape \cite{McMurry:2004aa} which has side lobes that decay exponentially fast.   The simulations in this article invoke a simple trapezoidal shape defined as
\begin{equation}
\label{eq: kappa2}
\kappa(t)=
\begin{cases}
1, &|t|\le \frac{1}{2}\\
-2|t|+2, &\frac{1}{2}\le|t|\le 1\\
0, & \text{else}
\end{cases}
\end{equation}
Taking the Fourier transform of this function gives the infinite-order kernel of interest:
\begin{equation}
\label{eq: K}
K(x)=\frac{2\left(\cos(x/2)-\cos(x)\right)}{\pi x^{2}}.
\end{equation}

We demonstrate the performance of using this infinite-order kernel for randomly right censored density and harzard function estimation in finite sample simulations.  Reproducible code for all of the simulations are provided as supplementary materials.  

\subsection{Normal Kernel vs Infinite-Order Kernel with Normal Data}

In this simulation we simply compare the performance of a normal kernel against the infinite-order kernel \eqref{eq: K} and remove the complicating issue of bandwidth selection.  Specifically, we determine the MSE performance for each estimator under their respective optimal bandwidth.   Lifetime and censoring data is simulated independently following a standard normal distribution thus yielding a censoring rate of 50\% on average.  The characteristic function of the standard normal distribution is $\phi(t)=\exp(-\frac{t^2}{2})$, which implies Assumption B is valid and the infinite-order kernel is asymptotically more efficient.  Estimates of the normal density at three points ($x$=0, 1, and 2) are considered along with two different sample sizes ($n$=50 and 500).  Results are provided for 999 realizations, which is sufficiently large to yield very small confidence intervals of the estimates.  The results of the simulation study (Table 1) shows improved MSE performance when using an infinite-order kernel, particularly with the larger sample size.  

\begin{table}[H]
\label{tab:sim1}
\caption{Comparison of the infinite-order kernel to the normal kernel with their respective optimal bandwidths.}
\centering
\begin{tabular}{ccccc}
\hline
 &  & \underline{$x=0$} & \underline{$x=1$} & \underline{$x=2$} \\
\multirow{2}{*}{$n=50$} & $\text{MSE}_{\text{infinite}}$$^{*}$ & $\textbf{3.96}_{.40}$ &  $\textbf{1.98}_{.70}$ &  $1.78_{.50}$\\
& $\text{MSE}_{\text{normal}}$$^{*}$ & $5.90_{.50}$ &  $3.93_{.90}$ &  $\textbf{1.33}_{.90}$  \\
\hline
\hline
 & & \underline{$x=0$} & \underline{$x=1$} & \underline{$x=2$} \\
\multirow{2}{*}{$n=500$} & $\text{MSE}_{\text{infinite}}$$^{*}$ & $\textbf{.54}_{.30}$ &  $\textbf{.28}_{.50}$ &  $\textbf{.47}_{.40}$\\
& $\text{MSE}_{\text{normal}}$$^{*}$ & $1.14_{.30}$ &  $.60_{.50}$ &  $.61_{.50}$  \\
\hline
\multicolumn{5}{l}{*\footnotesize{MSE values are multiplied by $10^3$ for easier comparison}}\\
\multicolumn{5}{l}{\footnotesize{and subscripted values correspond to the optimal bandwidth.}}\\
\end{tabular}
\end{table}

\subsection{Hazard Function Estimation With $\chi^2$ Data}

In this simulation we evaluate the performance of kernel density estimation on $\chi^2_\nu$ data using three different degrees of freedom: $\nu=7$, $\nu=11$, and $\nu=15$.  The characteristic function for the $\chi^2_\nu$ distribution is $(1-2it)^{-\nu/2}$, which implies assumption $A(r)$ holds for $r=2$, 4, and 6, respectively.  Again, two sample sizes ($n=50$ and $n=500$) are considered.  We first demonstrate the performance of the adaptive bandwidth selection algorithm discussed in Section \ref{sec:abc}.

\subsubsection{Bandwidth Selection Algorithm}
\begin{figure}[h]
\begin{center}
\includegraphics[width=.75\textwidth]{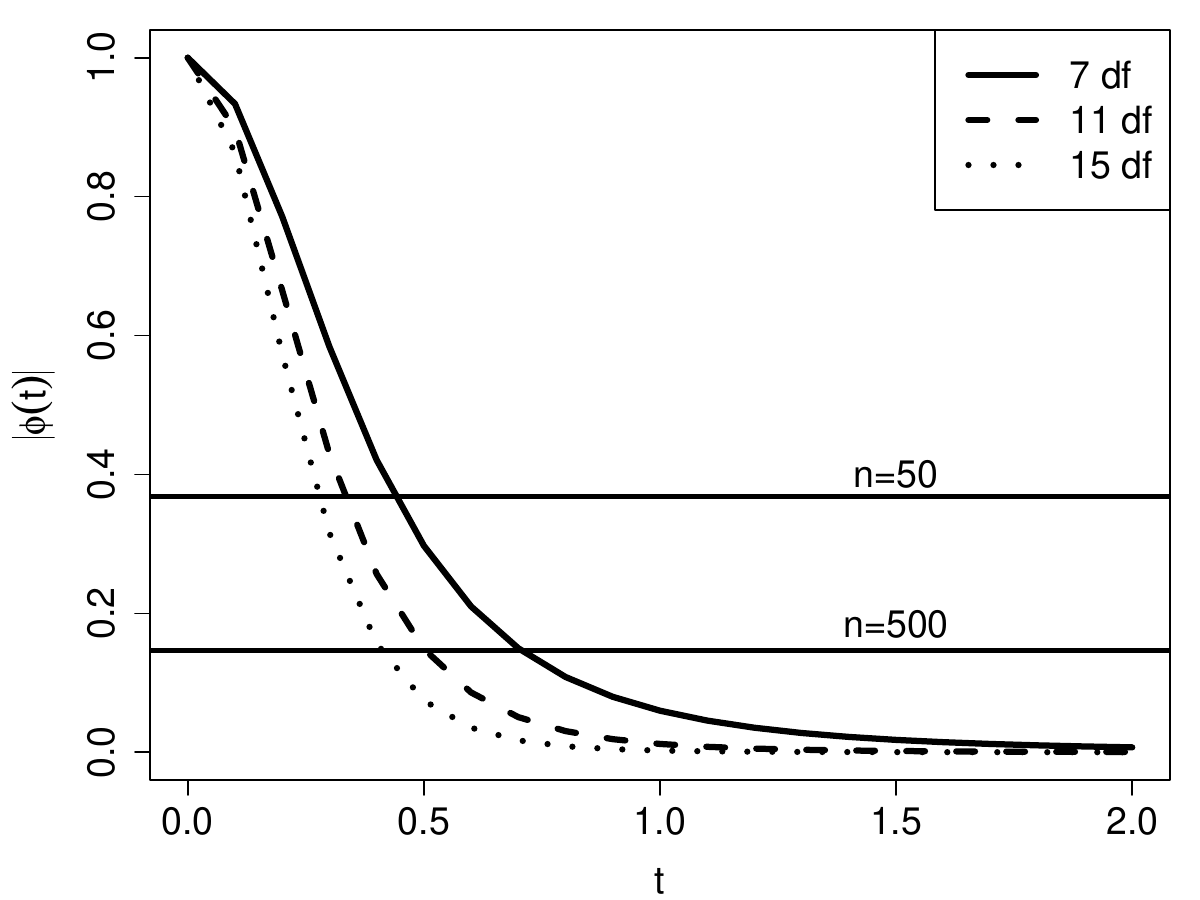}
\end{center}
\caption{The true characteristic function for each of the three $\chi^2$ densities with horizontal lines corresponding to the threshold in \protect\eqref{eq: alg}.}
\label{fig:phigraph}
\end{figure}

The true characteristic function for each of the three densities are graphed in Figure \ref{fig:phigraph}.  The two horizontal lines correspond to the thresholds given in Equation \eqref{eq: alg} for $C=2$ and $n=50$ and 500 respectively.  Following the bandwidth selection algorithm, we let $t^*$ be the value of $t$ corresponding to the point where $|\phi(t)|$ crosses the horizontal line, and then we set $h=1/(2t^*)$.   

\begin{figure}[h]
\begin{center}
\includegraphics[width=.49\textwidth]{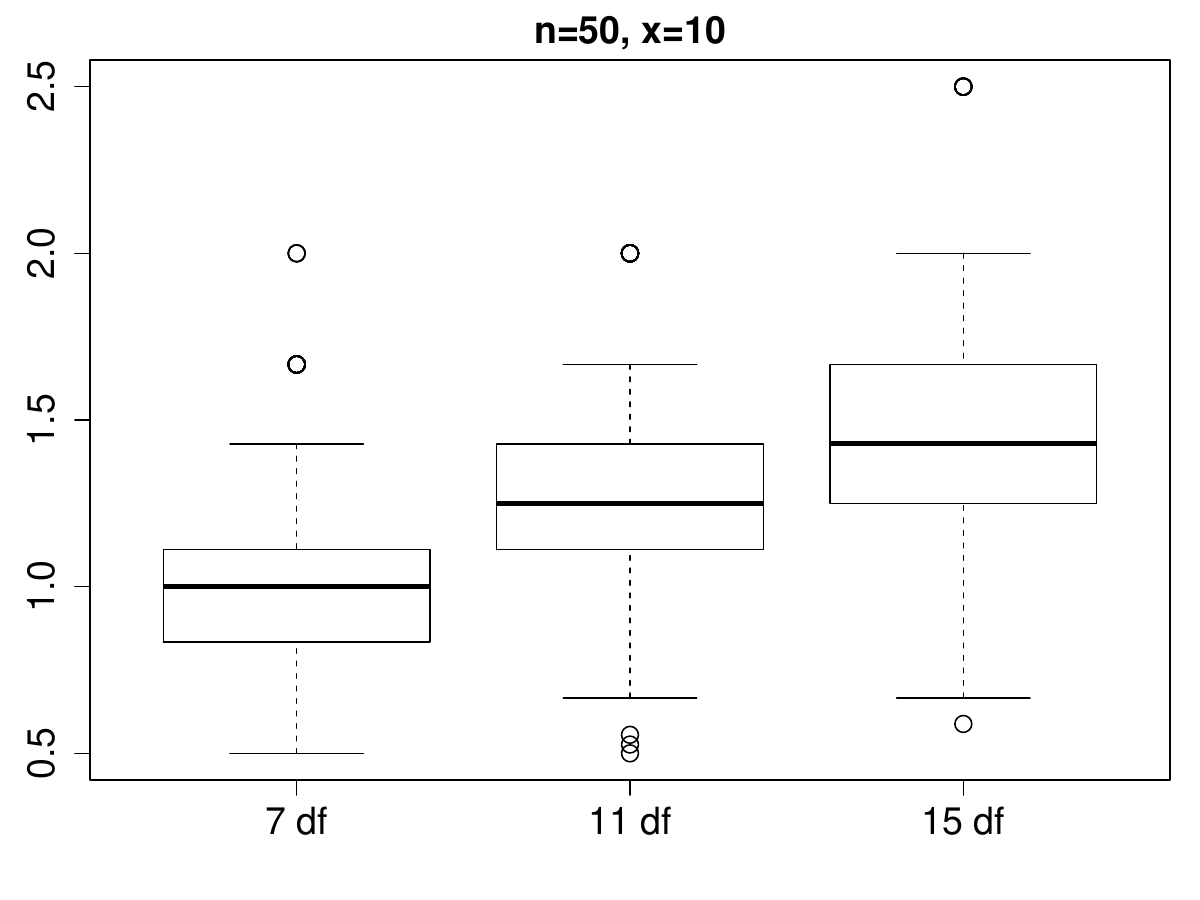}
\includegraphics[width=.49\textwidth]{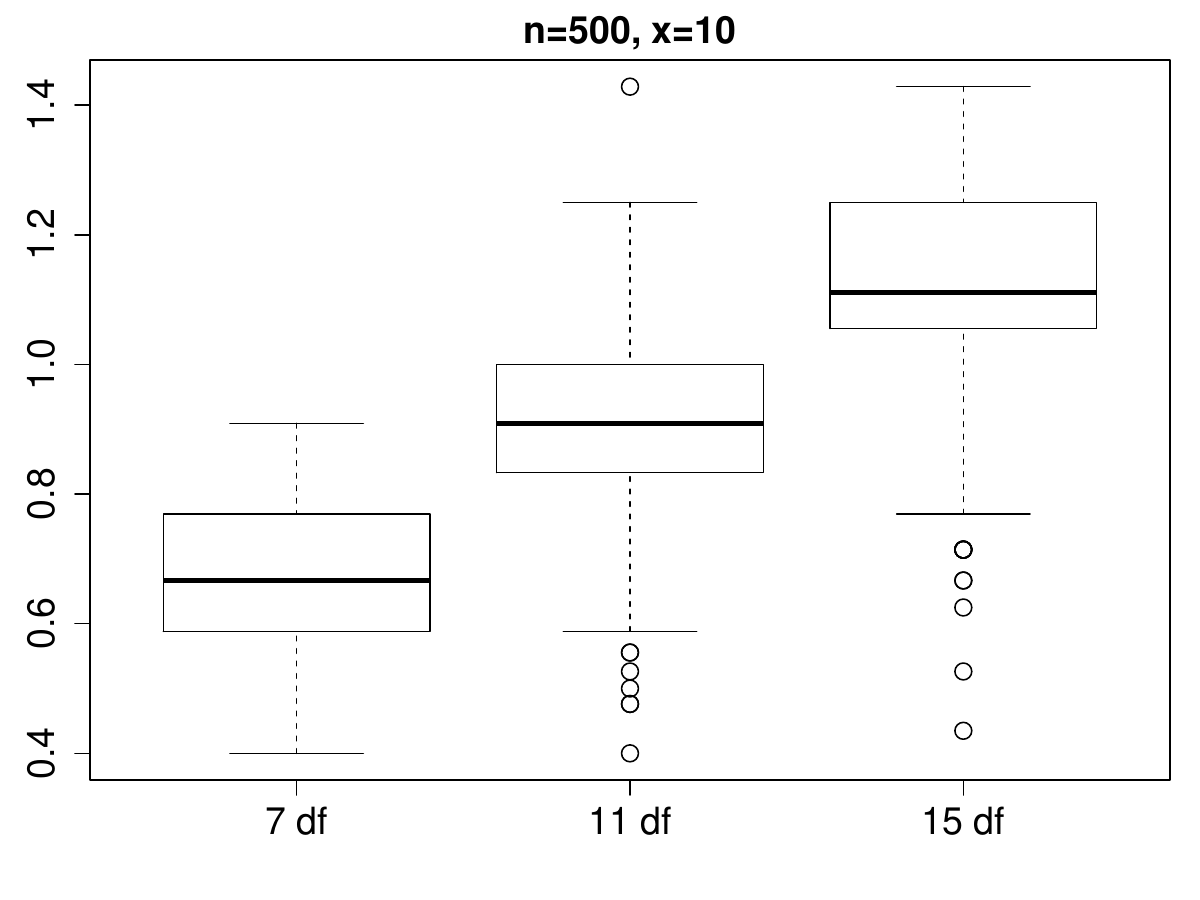}
\end{center}
\caption{The distribution of bandwidths from the bandwidth selection algorithm for estimating the density at $x=10$ when n=50 (left) and n=500 (right).}
\label{fig:actual}
\end{figure}

Figure \ref{fig:actual} shows the distribution of bandwidths for estimating the density at $x=10$ as determined by the bandwidth selection algorithm.  This fully automated procedure consistently identified the bandwidths in a narrow range, and its adaptive nature is observed as it produces increasingly larger bandwidths as the smoothness of the underlying density increases.  It also adapts to the sample size by producing smaller bandwidths with larger sample sizes.  

%

%

\subsubsection{Comparison of Hazard Estimators}

In many situations, particularly involving censored data, the support is known to lie in a half-line, or some compact interval, and
unaltered versions of kernel density estimators are not consistent near the boundary points.  However, a number of fixes for this boundary issue are available (see \cite{jones93} for a survey of several methods), and we adopt the simple reflection principle to resolve boundary problems in our estimator.  Specifically, when the density is known to have its support on $[0,\infty)$, we use the estimator $\hat{\hat{f}}(x)=\hat{f}(x)+\hat{f}(-x)$ to ensure consistency near the boundary point $x=0$; see \cite{schuster85} and \cite{silverman86} for discussions of this method with noncensored
data. 

In Table 2 we compare various estimators of the hazard function on the $\chi^2$ data.  The infinite-order kernel estimator of the hazard function is $\hat{f}(x)/\tilde{S}(x)$ where $\hat{f}(x)$ is the usual infinite-order density estimator and $\tilde{S}(x)$ is a smoothed Kaplan-Meier estimator (the R function \texttt{ksmooth} was applied to $\hat{S}$ to produce $\tilde{S}(x)$).  The other estimators considered are derived from the R packages \texttt{muhaz} and \texttt{survPresmooth}.  The \texttt{muhaz} estimator is based on the paper \cite{Muller:1994aa} with local bandwidth selection (denoted muhaz-l).  The \texttt{presmooth} estimator is based on the paper \cite{Tanner:1983aa} and uses the plug-in method for bandwidth selection \cite{Cao:2007aa}.

\begin{table}[H]
\label{tab:sim2}
\caption{Comparison of the infinite-order kernel to the normal kernel with their respective optimal bandwidths.}
\centering
\begin{tabular}{ccccc}
\hline
 &  & \underline{infinite} & \underline{muhaz-l} & \underline{presmooth} \\
\multirow{2}{*}{$n=50$} & $7$ df & $\textbf{2.20}_{.40}$ &  $\textbf{4.33}_{.65}$ &  $\textbf{2.36}_{.50}$\\
& $11$ df & $3.04_{.50}$ &  $4.39_{1.00}$ &  $3.37_{.85}$  \\
& $15$ df & $3.04_{.50}$ &  $4.39_{1.00}$ &  $3.37_{.85}$  \\
\hline
\hline
 &  & \underline{infinite} & \underline{muhaz-l} & \underline{presmooth} \\
\multirow{2}{*}{$n=500$} & $7$ df & $\textbf{2.20}_{.40}$ &  $\textbf{4.33}_{.65}$ &  $\textbf{2.36}_{.50}$\\
& $11$ df & $3.04_{.50}$ &  $4.39_{1.00}$ &  $3.37_{.85}$  \\
& $15$ df & $3.04_{.50}$ &  $4.39_{1.00}$ &  $3.37_{.85}$  \\
\hline
\multicolumn{5}{l}{*\footnotesize{MSE values are multiplied by $10^3$ for easier comparison}}\\
\multicolumn{5}{l}{\footnotesize{and subscripted values correspond to the optimal bandwidth.}}\\
\end{tabular}
\end{table}


\section{Breast carcinoma survival data}
\cite{Yakovlev2000} analyzed hazard functions on survival data involving 13,166 breast carcinoma patients identified through the Utah Cancer Registry.  In their study, a piecewise hazard function and a kernel based method was used to estimate the hazard functions of the individuals across different strata based on age and whether carcinoma was localized or not.  This dataset is re-analyzed with the proposed hazard function estimator along with the \texttt{muhaz} and \texttt{pehaz} estimators of \cite{Hess:2014aa,Muller:1994aa}.   

\begin{figure}[h]
\begin{center}
\includegraphics[width=.5\textwidth]{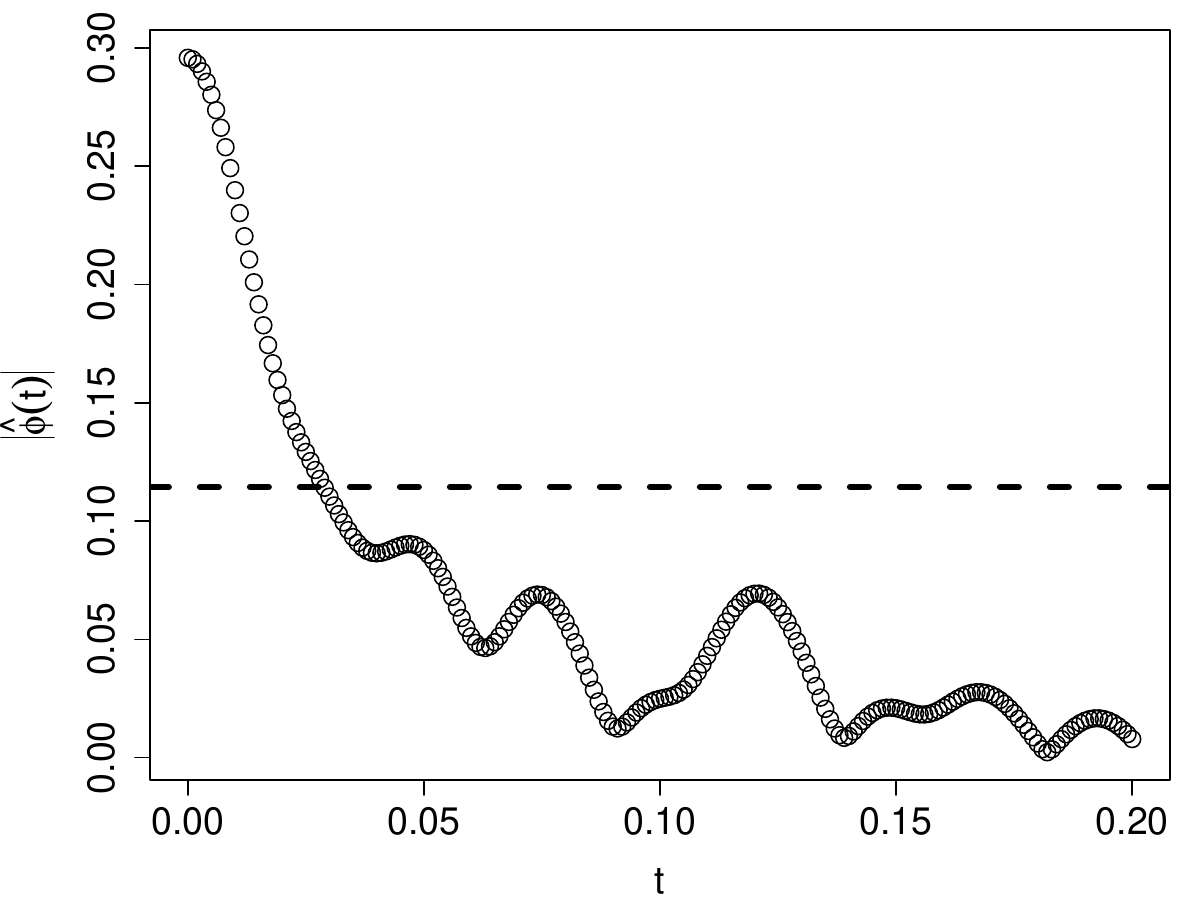}
\end{center}
\caption{$|\hat\phi(t)|$ for the breast cancer dataset (age$<$46, local) along with the threshold used to determine the bandwidth.}
\label{fig:phihatplot2}
\end{figure}

Figure \ref{fig:phihatplot2} shows a graph of $|\hat\phi(t)|$ along with the same threshold as used in the simulations.  We can observe $\hat\phi(t)$ smoothly decays toward zero in this real dataset.  This allows one to easily determine a reasonable range for the bandwidths to accompany the kernel density estimator.  

\begin{figure}[h]
\begin{center}
\includegraphics[width=.49\textwidth]{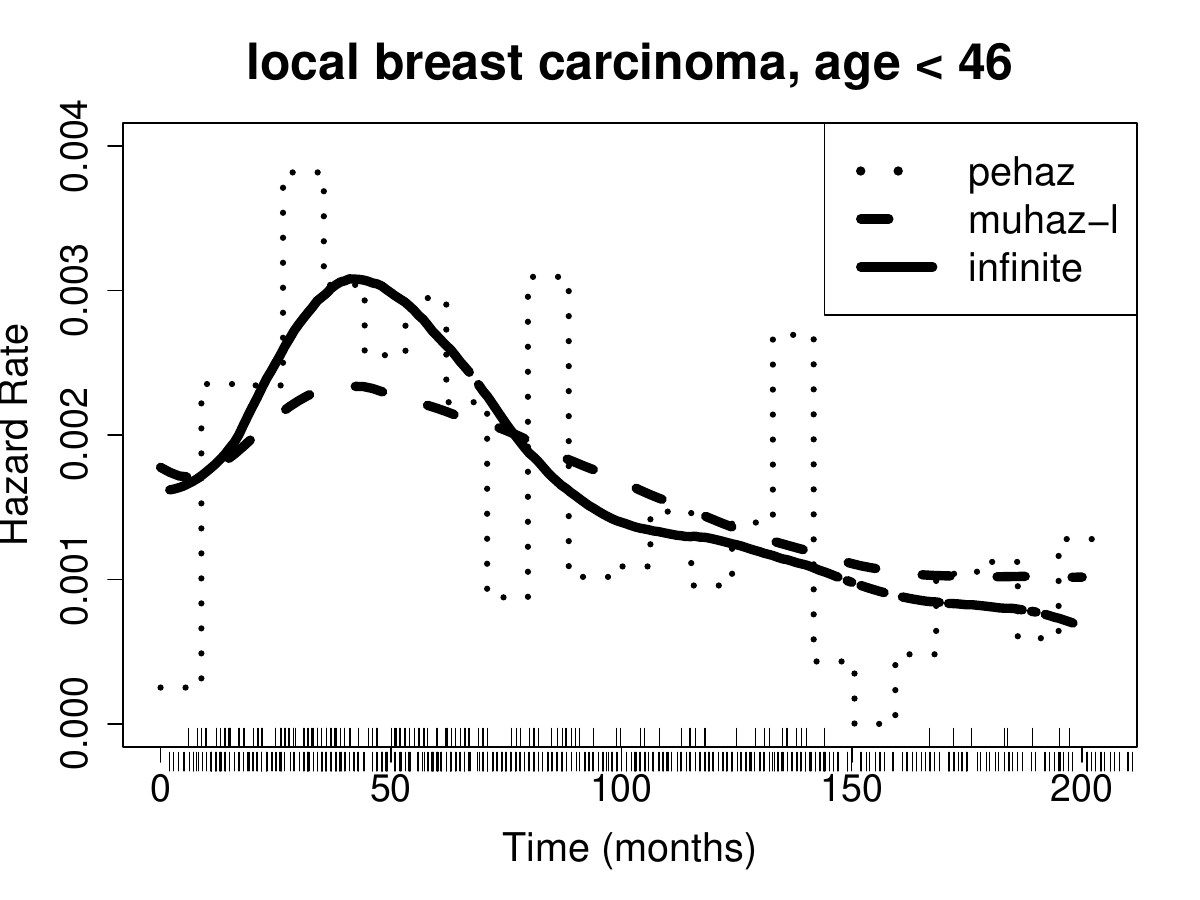}
\includegraphics[width=.49\textwidth]{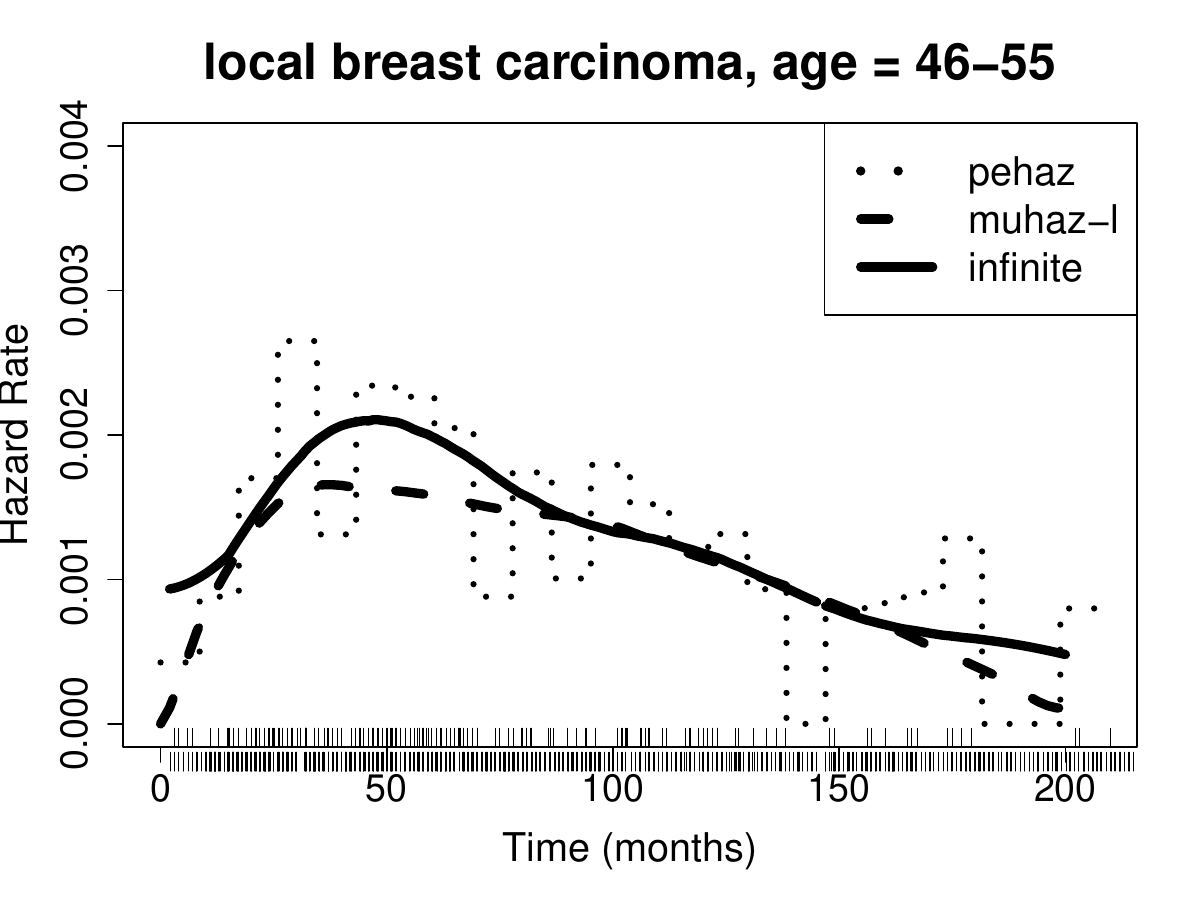}\\
\includegraphics[width=.49\textwidth]{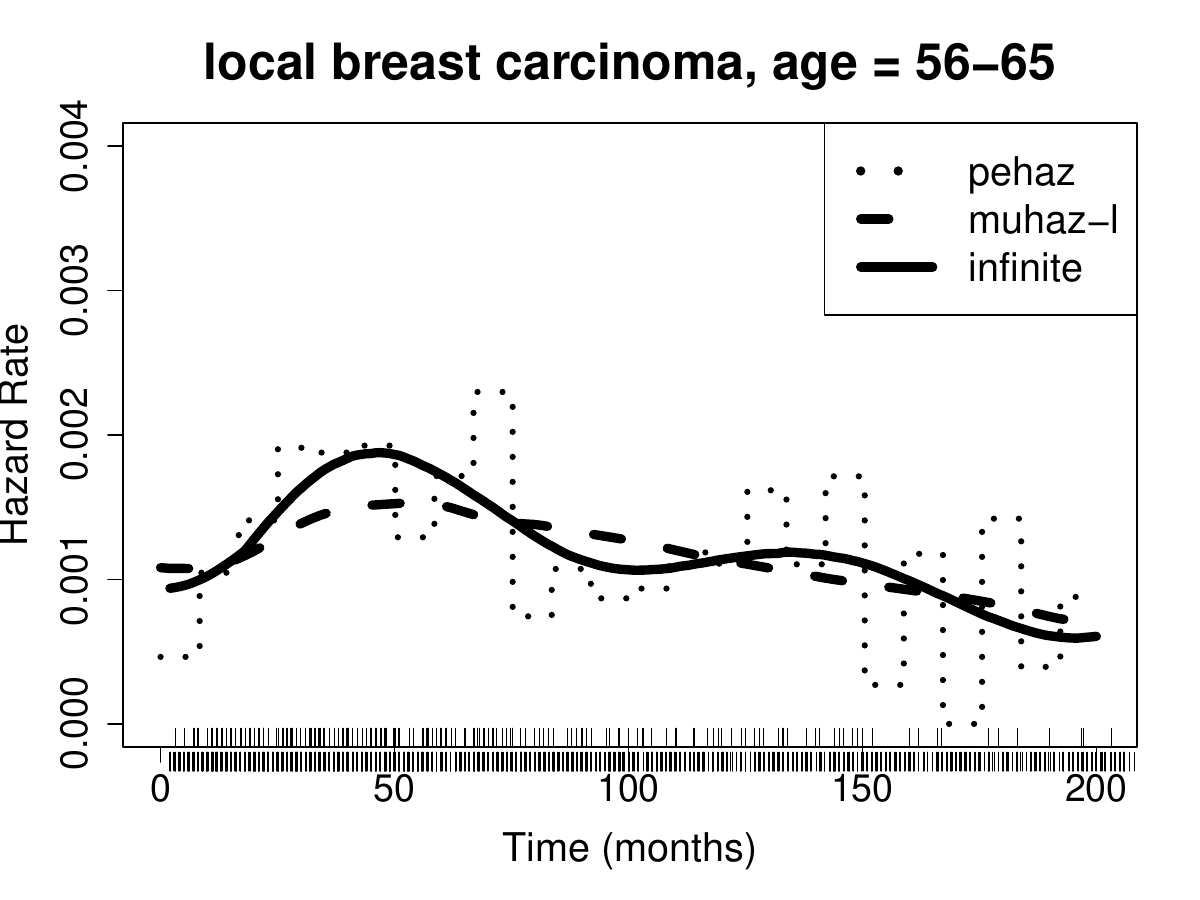}
\includegraphics[width=.49\textwidth]{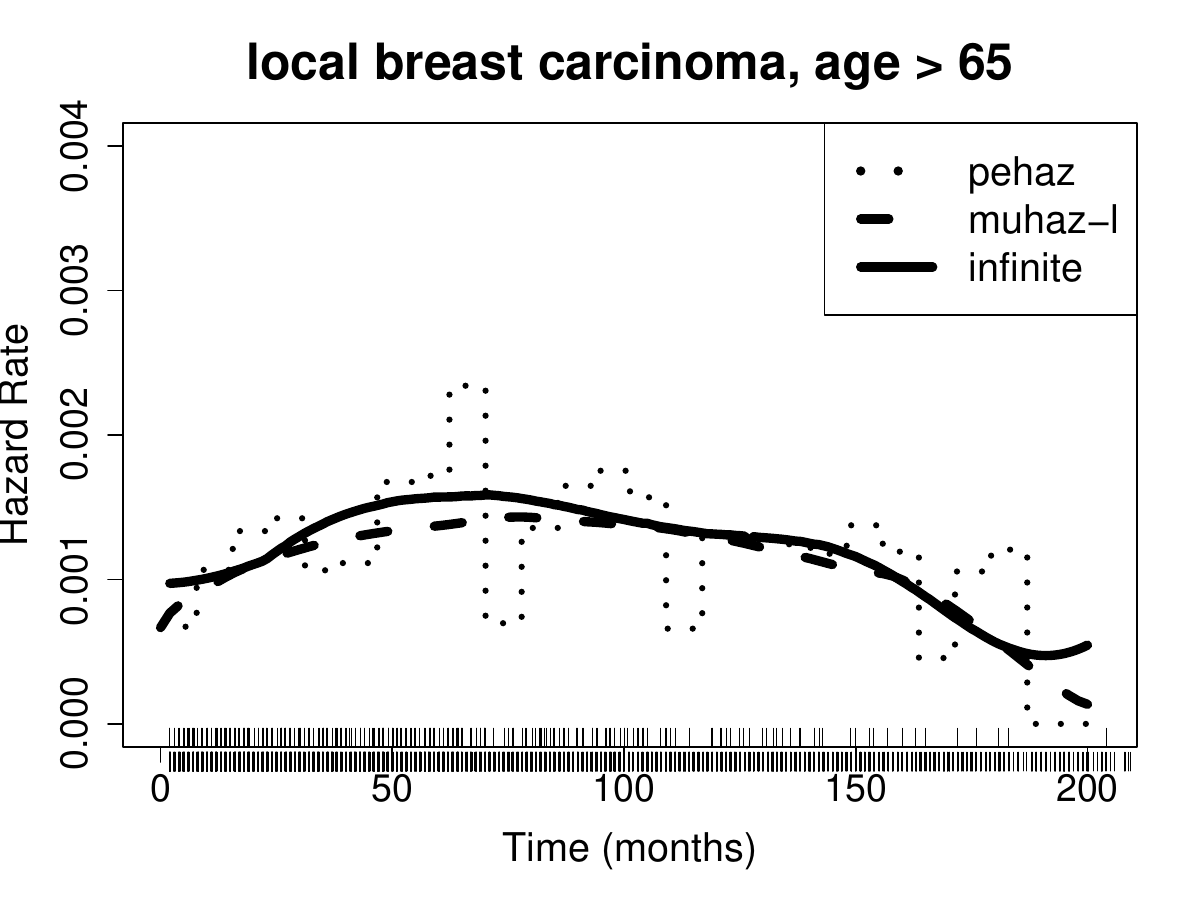}
\end{center}
\caption{The result of the three different hazard function estimators for localized breast cancer survival data for four different age ranges.}
\label{fig:local}
\end{figure}

\begin{figure}[h]
\begin{center}
\includegraphics[width=.49\textwidth]{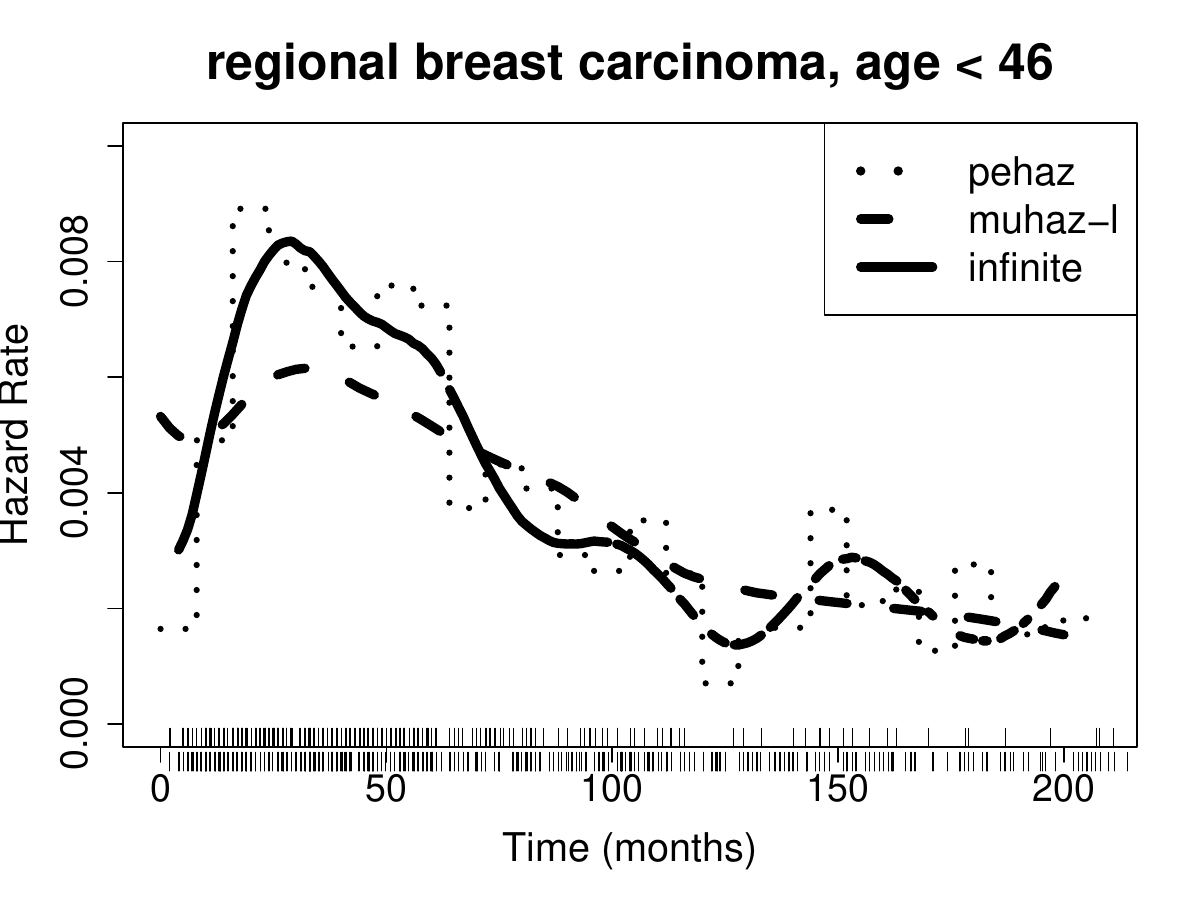}
\includegraphics[width=.49\textwidth]{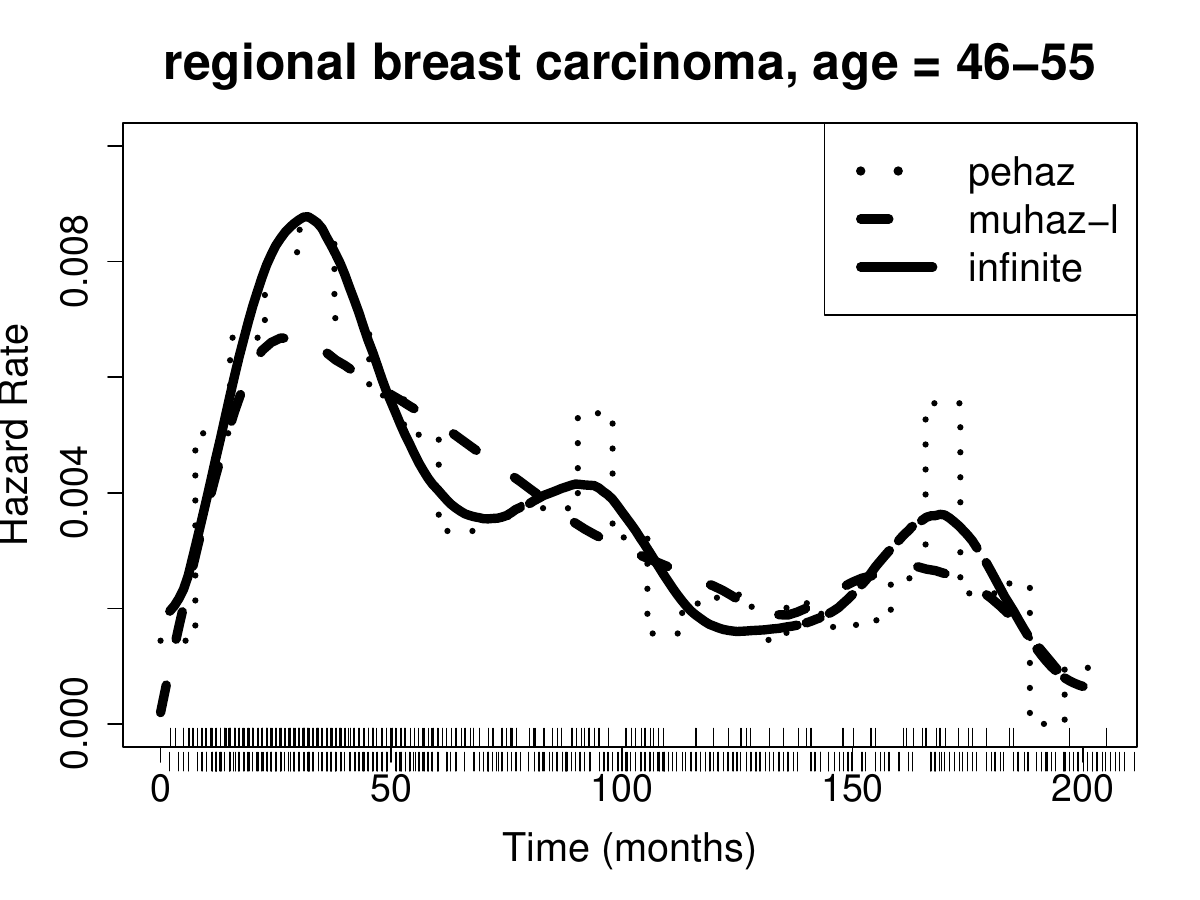}\\
\includegraphics[width=.49\textwidth]{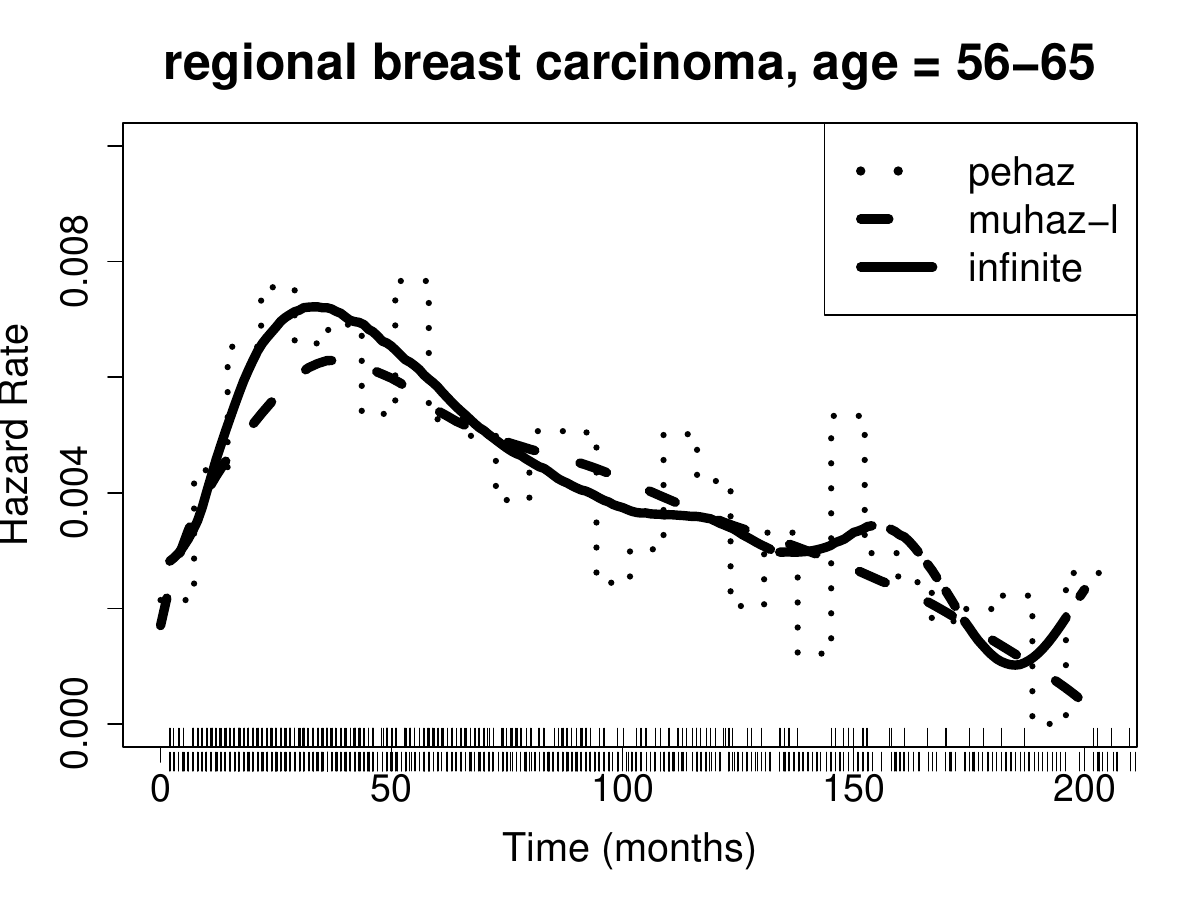}
\includegraphics[width=.49\textwidth]{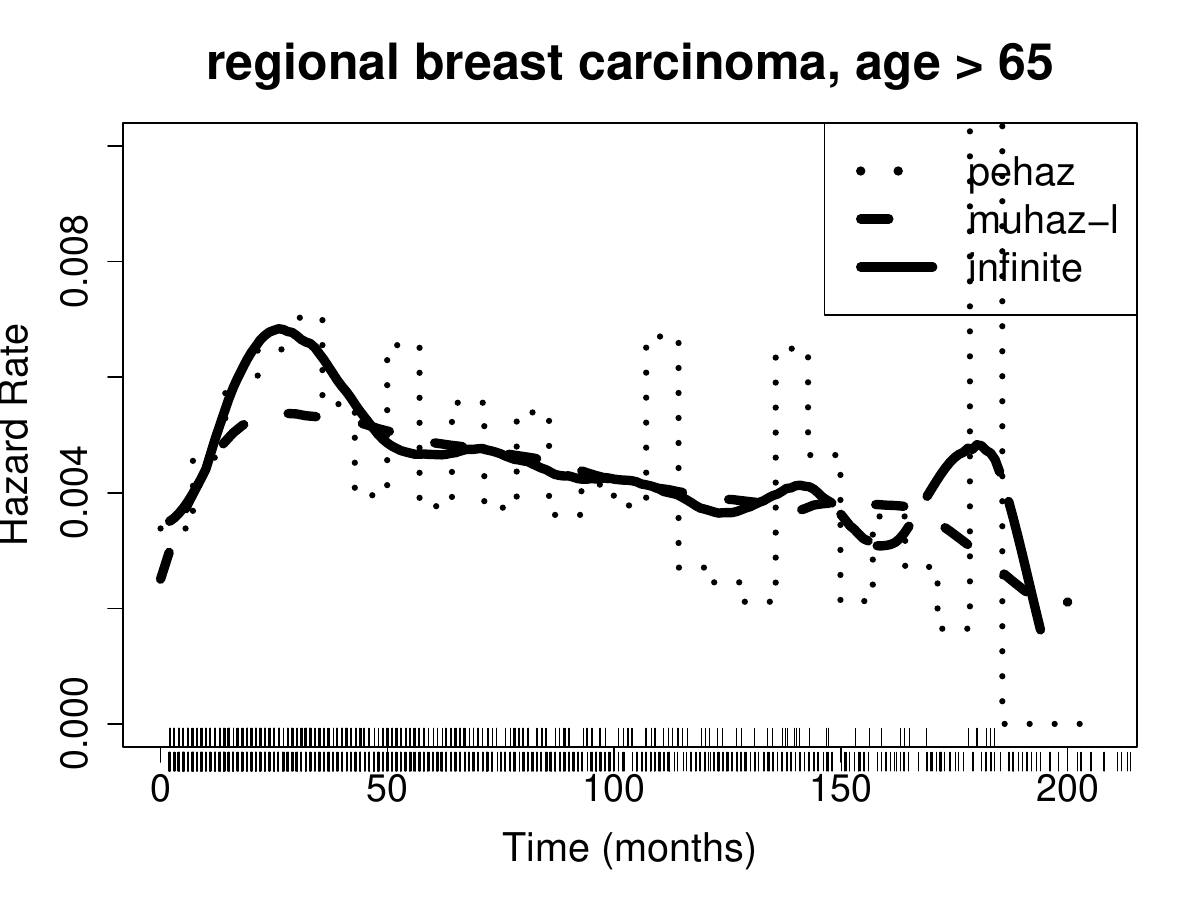}
\end{center}
\caption{The result of the three different hazard function estimators for non-local (regional) breast cancer survival data for four different age ranges.}
\label{fig:regional}
\end{figure}

Figures \ref{fig:local} and \ref{fig:regional} present the results of the different hazard function estimators on the breast cancer dataset.  It is consistently depicted among all of the estimators that as the severity of the disease increases, so does the hazard rate.  There is little difference in the estimated hazard rates for the different age groups.    The muhaz and infinite-order kernel estimator with adaptive bandwidth choice perform similarly on this dataset.  

\section{Conclusions}

The proposed infinite-order estimator, when used with its tailored bandwidth selection algorithm, produces a nearly $\sqrt{n}$-convergent
nonparametric estimator when the underlying density is sufficiently smooth, which corresponds to a rapidly decaying characteristic function.  Even in the least ideal situation of a slow decay of the characteristic function to zero (i.e., when the density is not very smooth), the estimator maintains the same performance as traditional kernel density estimators of censored data.  The same kernel was used throughout all of the simulations, so no
parameter estimation was involved in choosing the kernel, and the accompanying bandwidth selection algorithm requires very little computation to implement. Additionally, the proposed estimator is robust to sample size since no parameter estimation is involved and it can succeed in estimating the hazard function and density in small sample sizes where competing estimators may fail to produce an estimate.  Finally, the proposed estimator demonstrated reliable performance on the simulated data as well as on a actual datasets.

\appendix
\section{Technical Proofs}

\subsection{\enspace \textsc{Proof of Lemma \ref{lem:bias}.}}

 Theorem 2.1 in \cite{Zhou:1988aa} provides the following result: if $\theta(t)$ is a continuous nonnegative measurable function with $\E[\theta(X_{1})]<\infty$, then
\[
0\le \int_{-\infty}^{\infty} \theta(t)\,dF(t)-\E\left(\int_{-\infty}^{\infty} \theta(t)\,d\hat{F}(t)\right) \le \int_{-\infty}^{\infty} P(Z_{1}\le t)^{n} \theta(t)\,dF(t).
\]
By linearity of the integral and since $\theta(t)=\theta^{+}(t)-\theta^{-}(t)$ where $\theta^{+}(t)=\max(\theta(t),0)$ and $\theta^{-}(t)=\max(-\theta(t),0)$, we have the following result for general $\theta(t)$ 
\[
\left|\int_{-\infty}^{\infty} \theta(t)\,dF(t)-\E\left(\int_{-\infty}^{\infty} \theta(t)\,d\hat{F}(t)\right)\right| \le \int_{-\infty}^{\infty} P(Z_{1}\le t)^{n} (\theta^{+}(t)+\theta^{-}(t))\,dF(t).
\]
In particular, for $\theta(t)=e^{itx}=\cos(tx)+i\sin(tx)$, it follows that
\[
\begin{split}
\left|\text{bias}\left(\hat\phi(x)\right)\right|&=
\bigg|\int_{-\infty}^{\infty} \theta(t)\,dF(t)-\E\left(\int_{-\infty}^{\infty} \theta(t)\,d\hat{F}(t)\right)\bigg|\\
&\le \left|\int_{-\infty}^{\infty} \cos(tx)\,dF(t)-\E\left(\int_{-\infty}^{\infty} \cos(tx)\,d\hat{F}(t)\right)\right|\\
&\quad+\left|\int_{-\infty}^{\infty} \sin(tx)\,dF(t)-\E\left(\int_{-\infty}^{\infty} \sin(tx)\,d\hat{F}(t)\right)\right|\\
&\le 2\int_{-\infty}^{\infty} P(Z_{1}\le t)^{n}\,dF(t).\\
&\le 2 \max \left(\int_{-\infty}^{\infty}F(t)^{n}\,dF(t),\int_{-\infty}^{\infty}G(t)^{n}\,dF(t)\right)
\end{split}
\]
Note that
\[
\int_{-\infty}^{\infty}F(t)^{n}\,dF(t)=\frac{F(t)^{n+1}}{n+1}\biggl|_{-\infty}^{\infty}=O\left(\frac{1}{n}\right)
\]
From the assumptions of the lemma, we have $f(x)/g(x)\le M$ for some $M>0$, which gives
\[
\begin{split}
\int_{-\infty}^{\infty}G(t)^{n}\,dF(t)&=\int_{\{t:f(t)\not=0\}}G(t)^{n}f(t)\,dt\\
&=\int_{\{t:f(t)\not=0\}}G(t)^{n}\frac{f(t)}{g(t)}g(t)\,dt\\
&\le M\int_{\{t:f(t)\not=0\}}G(t)^n\,dG(t)\\
&=M\frac{G(t)^{n+1}}{n+1}\biggl|_{-\infty}^{\infty}\\
&=O\left(\frac{1}{n}\right)
\end{split}
\]

This establishes the bias of $\hat{\phi}(t)$ is $O(1/n)$ under the assumptions of Lemma \ref{lem:bias}.

\subsection{\textsc{Proof of Theorem \ref{thm: cen-bias}.}}

\begin{proof}
In order to evaluate the bias of $\hat{f}(x)$, we reformulate $\hat{f}(x)$ in terms of $\hat\phi(x)$ as follows
\begin{equation}
\label{eq: fhat2}
\begin{split}
  \hat{f}(x)
  &=\frac{1}{h}\sum_{j=1}^{n}s_j K\left(\frac{x-X_j}{h}\right)\\
  &=\frac{1}{h}\sum_{j=1}^n s_j \frac{h}{2\pi}\int_{-\infty}^\infty \kappa(t
  h)e^{-it(x-X_j)}\,dt\\
  &=\frac{1}{2\pi}\int_{-\infty}^\infty\left(\sum_{j=1}^ns_j
  e^{itX_j}\right)\kappa(t h)e^{-itx}\,dt\\
  &=\frac{1}{2\pi}\int_{-\infty}^\infty\hat{\phi}(t)\kappa(t
  h)e^{-itx}\,dt.
\end{split}
\end{equation}

From the representation in \eqref{eq: fhat2}, the expectation of $\hat{f}(x)$ is
\[
\begin{split}
\E[\hat{f}(x)]&=\frac{1}{2\pi}\int_{-\infty}^\infty\E\left[\hat{\phi}(t)\right]\kappa(th)e^{-itx}\,dt\\
&=\frac{1}{2\pi}\int_{-\infty}^{\infty} \phi(t)\kappa(th)e^{-itx}\,dt+O\left(\frac{1}{n}\right).
\end{split}
\]
Since $\phi(t)$ is the inverse Fourier transform of $f(x)$, $f(x)$
is therefore the Fourier transform of $\phi(t)$; i.e.
\begin{equation}
\label{eq: F(phi)}
f(x)=\frac{1}{2\pi}\int_{-\infty}^{\infty}\phi(t)e^{-itx}\,dt.
\end{equation}
Therefore the bias of $\hat{f}(x)$ is
\[
\bias(\hat{f}(x))=\E[\hat{f}(x)]-f(x)=\frac{1}{2\pi}\int_{-\infty}^\infty
(\kappa(th)-1)\phi(t)e^{-itx}\,dt + O\left(\frac{1}{n}\right).
\]
But since $\kappa(th)=1$ for $|t|\le 1/h$, we have
\[
\bias(\hat{f}(x))=\frac{1}{2\pi}\int_{|t|>1/h}(\kappa(th)-1)\phi(t)e^{-itx}\,dt + O\left(\frac{1}{n}\right).
\]
Since $|\kappa(t)|\le 1$ for all $t$, $|\kappa(th)-1|\le 2$ for all
$h$ and $t$.  We can then bound the bias by
\[
|\bias(\hat{f}(x))|\le \frac{2}{2\pi}\int_{|t|>1/h}|\phi(t)|\,dt + O\left(\frac{1}{n}\right).
\]

Under the assumption $\int |t|^r|\phi(t)|\,dt<\infty$ in (i), we
have
\[
\begin{split}
\int_{|t|>1/h}|\phi(t)|\,dt&=\int_{|t|>1/h}\frac{|t|^r|\phi(t)|}{|t|^r}\,dt\\
&=\le h^r\int_{|t|>1/h}|t|^r|\phi(t)|\,dt\\
&=o(h^r).
\end{split}
\]
If the bias is $o(h^r)+ O\left(\frac{1}{n}\right)$ and the variance is
$O\left(\frac{1}{nh}\right)$, then we wish to choose $h$ such that
$h^{2r}\sim\frac{1}{nh}$ which occurs if $h\sim an^{-\beta}$ with
$\beta=(2r+1)^{-1}$.  With this choice of $h$, we have
\[
\sup_{x\in\bR}
\left|\bias\left\{\hat{f}(x)\right\}\right|=o\left(n^{\frac{-r}{2r+1}}\right)\quad\text{and}\quad
\MSE\left\{\hat{f}(x)\right\}=O\left(n^{\frac{-2r}{2r+1}}\right).
\]
This proves part (i).

Under the assumption $|\phi(t)|\le De^{-d|t|}$ for some positive
constants $d$ and $D$, we have
\[
\begin{split}
\int_{|t|>1/h}|\phi(t)|\,dt&\le D\int_{|t|>1/h} e^{-d|t|}\,dt\\
&=\frac{D}{e^{d/h}}\int_{|t|>1/h}e^{d(1/h-|t|)}\,dt\\
&=O\left(e^{-d/h}\right)
\end{split}
\]
So the bias is $O(e^{-d/h})+ O\left(\frac{1}{n}\right)$, and by letting $h\sim1/(a\log n)$
gives a squared-bias of
\[
O\left(e^{\frac{-2d}{h}}\right)+ O\left(\frac{1}{n^{2}}\right)=O\left(e^{-2 d a\log
n}\right)+ O\left(\frac{1}{n^{2}}\right)=O\left(n^{-2da}\right)+ O\left(\frac{1}{n^{2}}\right)
\]
and a variance of
\[
O\left(\frac{1}{nh}\right)=O\left(\frac{a\log n}{n}\right).
\]
Therefore if $a>1/(2d)$, then
\[
\sup_{x\in\bR}
\left|\bias\left\{\hat{f}(x)\right\}\right|=O\left(\frac{1}{\sqrt{n}}\right)\quad\text{and}\quad
\MSE\left\{\hat{f}(x)\right\}=O\left(\frac{\log n}{n}\right)
\]
This proves part (ii).

Under the assumption $\phi(t)=0$ when $|t|\ge b$, we have
\[
\int_{|t|>1/h}|\phi(t)|\,dt=0
\]
when $h\le 1/b$.  So by letting $h\le 1/b$, we have
\[
\sup_{x\in\bR}
\left|\bias\left\{\hat{f}(x)\right\}\right|=O\left(\frac{1}{n}\right)\quad\text{and}\quad
\MSE\left\{\hat{f}(x)\right\}=O\left(\frac{1}{n}\right)
\]
which completes the proof of the theorem.
\end{proof}

\subsection*{A.2\enspace \textsc{Proof of Theorem \ref{thm: cen-bias2}.}}
\begin{proof} 
By taking the $p^{\text{th}}$ derivative on both sides of the
identity \eqref{eq: identity}, we have
\[
\frac{1}{h^{p+1}}K^{(p)}\left(\frac{x}{h}\right)=\frac{1}{2\pi}\int_{-\infty}^\infty
(-i t)^p \kappa(t h)e^{-itx}\,dt.
\]
By taking the $p^{\text{th}}$ derivative on both sides of the
identity \eqref{eq: F(phi)}, we have
\[
f^{(p)}(x)=\frac{1}{2\pi}\int_{-\infty}^\infty(-it)^p
\hat{\phi}(t)\kappa(th)e^{-itx}\,dt.
\]
Following the steps in \eqref{eq: fhat2}, we have
\[
\hat{f}_p(x)=\frac{1}{2\pi}\int_{-\infty}^\infty\hat{\phi}(t)\kappa(th)e^{-itx}\,dt.
\]
and we can now compute the bias of $\hat{f}_p(x)$ to be
\[
\bias(\hat{f}_p(x))=\frac{1}{2\pi}\int_{-\infty}^\infty
(-it)^p(\kappa(th)-1)\phi(t)e^{-itx}\,dt + O\left(\frac{1}{n}\right).
\]
Proceeding as in the proof of Theorem \ref{thm: cen-bias}, this bias
is bounded as
\[
|\bias(\hat{f}_p(x))|\le
\frac{2}{2\pi}\int_{|t|>1/h}|t|^p|\phi(t)|\,dt + O\left(\frac{1}{n}\right).
\]
Under assumption $A(r+p)$, we have
\[
\begin{split}
\int_{|t|>1/h}|t|^p|\phi(t)|\,dt&=\int_{|t|>1/h}\frac{|t|^{r+p}\phi(t)}{|t|^{r}}\,dt\\
&\le h^r\int_{|t|>1/h}|t|^{r+p}|\phi(t)|\,dt\\
&=o(h^r).
\end{split}
\]
If the bias is $o(h^r) + O\left(\frac{1}{n}\right)$ and the variance is
$O\left(\frac{1}{nh^{p+1}}\right)$, then we wish to choose $h$ such
that $h^{2r}\sim\frac{1}{nh^{p+1}}$ which occurs if $h\sim
an^{-\beta}$ with $\beta=(2r+p+1)^{-1}$.  With this choice of $h$,
we have
\[
\sup_{x\in\bR}
\left|\bias\left\{\hat{f}_p(x)\right\}\right|=o\left(n^{\frac{-r}{2r+p+1}}\right)\quad\text{and}\quad
\MSE\left\{\hat{f}_p(x)\right\}=O\left(n^{\frac{-2r}{2r+p+1}}\right).
\]
Under assumption $B$,
\[
\begin{split}
\int_{|t|>1/h}|t|^p|\phi(t)|\,dt&\le D\int_{|t|>1/h}|t|^p e^{-d|t|}\,dt\\
&=\frac{D}{e^{d/h}}\int_{|t|>1/h}|t|^p e^{d(1/h-|t|)}\,dt\\
&=O\left(e^{-d/h}\right).
\end{split}
\]
Under assumption $C$,
\[
\int_{|t|>1/h}|t|^p|\phi(t)|\,dt=0
\]
when $h\le 1/b$.  Finally, the bias and MSE results for parts (ii)
and (iii) now follow along the same lines as Theorem \ref{thm:
cen-bias}.
\end{proof}

\subsection*{A.3\enspace \textsc{Proof of Theorem \ref{thm: cen-band}.}}
\begin{proof}
The proof follows the proof of Theorem \ref{thm:
cen-band} in \cite{Berg:2009aa} with
  little modification.
\end{proof}

\subsection*{A.4\enspace \textsc{Proof of Theorem \ref{thm: cen-pilot}.}}
\begin{proof}
Parts (ii) and (iii) follow from Theorems 1 and 2 and the
$\delta$-method.  The convergence of $\hat{h}_{\text{M}}$ in part
(i) is dictated by the slowly converging $\widehat{f''}(x)$.
However, the convergence rate of $\hat{h}_{\text{M}}$ is unhampered
by the convergence rate of $\hat{h}$; for instance, if $h$ is
replaced with the random quantity $h(1+o_p(1))$ (refer to the proof
of Lemma 2 in \cite{buhlmann96}) then Theorem \ref{thm: cen-bias} is
still valid.  If $|\phi(t)|\sim A|t|^{-d}$, then by Theorem \ref{thm:
cen-band},
\[
\hat{h}\stackrel{P}{\sim} \tilde{A}\left(\frac{\log
n}{n}\right)^{\frac{1}{2d}}
\]
From Theorem \ref{thm: derivative}, part (i), if
\begin{equation}
\label{eq: condition} \int_{-\infty}^\infty
|t|^{r+2}\,|\phi(t)|<\infty,
\end{equation}
then the bias of $\widehat{f''}(x)$ is $o(h^r)$.  In order for
\eqref{eq: condition} to be satisfied, $r$ must be less than $d-3$,
so we let $r=\lceil d-4\rceil$.  Therefore the bias of
$\widehat{f''}(x)$ (which dominates the MSE of $\widehat{f''}(x)$)
is
\[
o\left( \frac{\log n}{n}\right)^{\frac{\lceil d-4\rceil}{2d}},
\]
and coupled with the $\delta$-method, part (i) of Theorem 4 is now proved.
\end{proof}

\begin{acknowledgements}
We appreciate the efforts of Professors Alex Tsodikov and Richard Keeber who helped us to obtain the Utah Cancer Registry data that was analyzed in this manuscript. 
\end{acknowledgements}


\bibliographystyle{unsrt}
\bibliography{censored}

\end{document}